\documentclass[11pt]{article}
\usepackage{amsfonts}
\usepackage{mathrsfs}
\usepackage{amsmath}
\usepackage{amssymb}
\usepackage{graphicx}
\usepackage{amsthm,latexsym}
\renewcommand{\paragraph}{\roman{paragraph}}
\newtheorem{theorem}{\scshape \mdseries  Theorem}[section]
\newtheorem{lemma}[theorem]{\scshape \mdseries  Lemma}
\newtheorem{coro}[theorem]{\scshape \mdseries  Corollary}
\newtheorem{conj}[theorem]{\scshape \mdseries  Conjecture}

\begin{document}

\title{\sf The signature of line graphs and power trees
 \thanks{Supported by National Natural Science Foundation of China (11071002, 11371028),
 Program for New Century Excellent Talents in University (NCET-10-0001),
 Key Project of Chinese Ministry of Education (210091),
Specialized Research Fund for the Doctoral Program of Higher Education (20103401110002),
Scientific Research Fund for Fostering Distinguished Young Scholars of Anhui University(KJJQ1001).}
}
\author{Long Wang,  Yi-Zheng Fan\thanks{Corresponding author. E-mail address: fanyz@ahu.edu.cn (Y.-Z. Fan), wanglongxuzhou@126.com (L. Wang)}\\
{\small  \it  School of Mathematics Sciences, Anhui University, Hefei 230601, P.R. China} \\
   }
\date{}
\maketitle

\noindent {\bf Abstract:}\ \
Let $G$ be a graph and let $A(G)$ be the adjacency matrix of $G$.
The signature $s(G)$ of $G$ is the difference between the positive inertia index and the negative inertia index of $A(G)$.
Ma et al. [Positive and negative inertia index of a graph,  Linear Algebra and its Applications 438(2013)331-341] conjectured that
$-c_3(G)\leq s(G)\leq c_5(G),$
where $c_3(G)$ and $c_5(G)$ respectively denote the number of cycles in $G$ which have length $4k+3$ and $4k+5$ for some integers $k \ge 0$,
and proved the conjecture holds for trees, unicyclic or bicyclic graphs.

It is known that $s(G)=0$ if $G$ is bipartite, and the signature is closely related to the odd cycles or nonbipartiteness of a graph from the existed results.
In this paper we show that the conjecture holds for the line graph and power trees.



\section{Introduction}
Throughout this paper we consider only simple graphs.
The {\it adjacency matrix} $A(G)=[a_{ij}]$ of a graph $G$ with vertex set
$V(G)=\{v_1,v_2,\ldots, v_n\}$ and edge set $E(G)$ is defined to be a
symmetric matrix of order $n$ such that $a_{ij}=1$ if $v_i$ is
adjacent to $v_j$, and $a_{ij}=0$ otherwise.
The {\it positive inertia index} $p(G)$, the {\it negative inertia index} $n(G)$ and the {\it nullity} $\eta(G)$ of $G$ are respectively defined to be
the number of positive eigenvalues, negative eigenvalues and zero eigenvalues of $A(G)$.
The {\it rank} of $G$, written as $r(G)$, is defined to be the rank of $A(G)$.
The {\it signature} of $G$, denoted by $s(G)$, is defined to be the difference $p(G)-n(G)$.
Obviously,  $p(G)+n(G)+ \eta(G)=|V(G)|$, $p(G) + n(G)=r(G)$ and $p(G)-n(G)=s(G)$.

Motivated by the discovery that the nullity of a graph is related to the stability
of the molecular represented by the graph \cite{atk} and the open problem of characterizing all singular graphs posed by Collatz \cite{col},
 many authors discuss the nullity of a graph and obtain a lot of interesting results.
  Here we particularly mention the results involved with the nullity of line graphs.
Sciriha \cite{sci} proved that all trees whose line graph is singular must have an even order.
Gutman and Sciriha \cite{gut} showed that the  nullity of the line graph of a tree is at most one.
Li et al. \cite{li} proved that the nullity of the line graph of a unicyclic graph with depth one is at most two.
Gong et al. \cite{gon} improved the above results as: the nullity of the line graph of a connected graph with $k$ induced cycles is at most $k+1$.

Recently some authors discuss a more general problem, that is, describing the positive or negative inertia index of graphs or weighted graphs, especially of trees or their line graphs, unicyclic or bicyclic graphs; see Ma et al. \cite{ma}, Li et al. \cite{lisc} and Yu et al. \cite{yu,yugh}.
In the paper \cite{ma} the authors posed a conjecture as follows, and proved the conjecture holds for trees, unicyclic or bicyclic graphs.
\begin{conj}\label{conj} \emph{\cite{ma}}
The inequality  $-c_3(G)\leq s(G)\leq c_5(G)$ possibly holds for any simple graph $G$,
where $c_3(G)$ and $c_5(G)$ denote respectively  the number of cycles having length $4k+3$ (or length $3$ modulo $4$) and
 the number of cycles having length $4k+5$ for some integers $k \ge 0$  (or length $1$ modulo $4$).
\end{conj}

\begin{theorem}\label{tree}
\emph{\cite{ma}} Let $G$ be  a tree, or a unicyclic graph, or a bicyclic graph. Then $-c_3(G)\leq s(G)\leq c_5(G)$.
\end{theorem}

\noindent A  weaker result was also given by Ma et al. \cite{ma} that $|s(G)| \le c_1(G)$ for any graph $G$, 
where $c_1(G)$ denotes the number of odd cycles of $G$, or $c_1(G)=c_3(G)+c_5(G)$.

When $G$ is bipartite, surely $s(G)=0$ and the conjecture holds in this case.
So, from Theorem \ref{tree} or Conjecture \ref{conj} (if it was true), 
we find that the signature is closely related to the odd cycles or nonbipartiteness of a graph.
In this paper we prove that the conjecture holds for the line graphs and power trees.

\section{Preliminaries}
We first introduce some notations.
Let $G$ be a graph and let $W  \subseteq V(G)$.
Denote by $G - W$ the subgraph of $G$ obtained by deleting the vertices in $W$ together with all edges incident to them.
If $G_1$ is a subgraph of $G$, we sometimes write $G-G_1$ instead of $G-V(G_1)$.
In particular, if $W=\{x\}$, we simply write $G-W$ as $G-x$.
If $G_1$ is an induced subgraph of $G$ and $x$ is a vertex of $G$ outside $ G_1$,
denote by $G_1+x$ the subgraph of $G$ induced by the the vertices of $G_1$ and $x$.

\begin{lemma} {\em \cite{ma}} \label{sper}
Let $G$ be a graph containing path with four vertices of degree $2$ as shown in Fig. 2.1. Let $H$ be the
graph obtained from $G$ by replacing this path with an edge. Then $p(G) = p(H)+2$, $n(G) = n(H)+2$,
$\eta(G) = \eta(H)$, and hence $s(G)=s(H)$.
\end{lemma}

\begin{center}
\vspace{3mm}
\includegraphics[scale=.6]{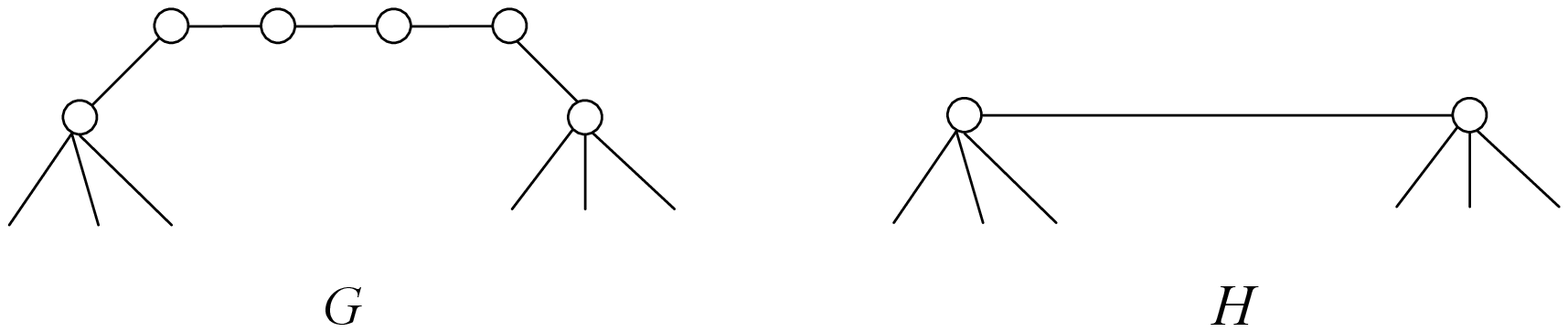} \label{f21}

{\small Fig. 2.1. The graphs $G$ and $H$ in Lemma \ref{sper}}
\end{center}

\begin{lemma}{\em \cite{li}} \label{luni}
Let $C_{n_1,n_2,\ldots,n_t}$ be the graph obtained from a cycle $C_t$ by attaching $n_i$ pendent edges to each vertex $v_i$ of $C_t$, where $n_i \ge 0$.
Let $G$ be the line graph of $C_{n_1,n_2,\ldots,n_t}$, and let $m=|\{i|n_i >0\}|$.
Then the following results hold, where a zero chain of finite integer sequence is defined to be a zero subsequence whose (cyclic) predecessor and successor are both nonzero, and the length of the zero chain is defined to be the number of integers in that zero subsequence.

{\em (1)} $\eta(G)=2$ if and only if $m=0$ and $t \equiv 0 \mod 4$.

{\em (2)} $\eta(G)=1$ if and only if $m \ge 1$ and either $n_i \in \{0,1\}$ for $i=1,2,\ldots,t$, the length of any zero chain of $(n_1,n_2,\ldots,n_t)$ is even,
and $t+m \equiv 0 \mod 4$; or $t \equiv 0 \mod 4$ and one of $n_1=n_3=\cdots=n_{t-1}=0$ and $n_2=n_4=\cdots=n_{t}=0$ must hold.

{\em (3)} $\eta(G)=0$ otherwise.
\end{lemma}

\begin{lemma}\label{cut vertex rank}\emph{\cite{gon}}  Let $x$ be a cut vertex of a graph $G$ and $G_{1}$ be a component of $G-x$.
If $r(G_1+x)=r(G_1)+2$, then $r(G)=r(G-x)+2$. If $r(G_1+x)=r(G_{1})$, then $r(G)=r(G_1)+r(G-G_1)$.\end{lemma}

\begin{lemma}\label{vertex-deleting sign 2}  Let $G$ be a graph and let $x$ be a vertex of $G$.
Then  $ |s(G)-s(G-x)|\leq 1$. In particular, if $r(G-x)=r(G)$ or $r(G-x)=r(G)-2$, then $s(G-x)=s(G)$.
\end{lemma}

\noindent
{\bf Proof.} By  the eigenvalues interlacing property of real symmetric matrices (or see \cite{cve}), we have $p(G)-1 \leq p(G-x)\leq p(G)$, and $n(G)-1 \leq n(G-x)\leq n(G)$, which yields the required results immediately. \hfill $\blacksquare$

\begin{coro}\label{induced}  Let $H$ be an induced subgraph of  a graph $G$. If $r(H)=r(G)$, then $s(H)=s(G)$.
\end{coro}

\noindent
{\bf Proof.} Note that $H$ can be viewed as one  obtained from $G$ by sequently deleting some of its vertices.
The condition $r(H)=r(G)$ implies that in each step the rank, and hence the signature of the resulting graph keeps invariant by Lemma \ref{vertex-deleting sign 2},
which yields $s(H)=s(G)$.  \hfill $\blacksquare$

\begin{coro}\label{cut vertex sign} Let $x$ be a cut vertex of a graph $G$ and let $G_1$ be a component of $G -x$.
If $r(G_1+ x)=r(G_1)+2$, then $s(G)=s(G-x)$.
\end{coro}

\noindent
{\bf Proof.}  If $r(G_1+ x)=r(G_1)+2$, then $r(G)=r(G-x)+2$ by Lemma \ref{cut vertex rank}, and hence $s(G)=s(G-x)$ by Lemma \ref{vertex-deleting sign 2}.
\hfill $\blacksquare$

\begin{lemma}\label{cut vertex sign 2}
Let $x$ be a cut vertex of a graph $G $ and let $G_1, G_2, \ldots, G_k$  be all components of $G -x$.
If $r(G_1) =r(G_1+x)$, then $s(G)=s(G_1)+s(G-G_1)$.
In particular, if $r(G_i) =r(G_i+x)$ for all $i$, then $s(G)=s(G-x)$.
\end{lemma}

\noindent{\bf Proof.}  Let $\Gamma=\cup_{i=2}^k G_k$. Write the adjacency matrix of $G$ as follows,
$$A(G)=\left[\begin{array}{ccc}
A(G_1)&\alpha&0\\
\alpha^T&0&\beta\\
0&\beta^T&A(\Gamma)
\end{array}\right],$$
where the middle $0$ corresponds to the cut vertex $x$.
 As $r(G_1)=r(G_1+x)$, the matrix equation $A(G_1)X=\alpha$ has a solution, say $\xi$, such that $\alpha^T\xi=0$.
 Now, take $Q$ as the following matrix with the same partition as $A(G)$,
 $$Q=\left[\begin{array}{ccc}
 I&-\xi&0\\
 0&1&0\\
0&0&I
\end{array}\right],$$
Then
$$Q^TA(G)Q=\left[\begin{array}{ccc}A(G_1)&0&0\\0&0& \beta\\
0&\beta^T&A(\Gamma)\end{array}\right].$$
So we have $s(G)=s(G_1)+s(G-G_1)$.

If $r(G_i) =r(G_i+x)$ for all $i$, by induction on the number of components of $G-x$, we have
$s(G)=\sum_{i=1}^{k-1} s(G_i)+s(G_k+x)$.
The result follows as $s(G_k+x)=s(G_k)$ by Lemma \ref{vertex-deleting sign 2}.\hfill $\blacksquare$

\begin{lemma}\label{cut vertex sign 3} Let $x$ be a cut vertex of a graph $G$ and $G_1, G_2, \ldots, G_k$  be all components of $G -x$.
If $s(G) = s(G-x)+1$, then $s(G_l+x)=s(G_l)+1$ for some $l$, and $s(G)=s(G_{l}+x)+\sum_{j\neq l}s(G_{j})$.
\end{lemma}

\noindent{\bf Proof.}
Note that $r(G_i+x) \le r(G_i)+2$ for each $i$.
If $r(G_i+x)=r(G_i)+2$ for some $i$ or $r(G_i+x)=r(G_i)$ for all $i$'s, then $s(G)=s(G-x)$ by Corollary \ref{cut vertex sign} or Lemma \ref{cut vertex sign 2}; a contradiction.
So $r(G_i+x) \le r(G_i)+1$ for all $i$'s, with equality for at least one $i$.

Write the adjacency matrix of $G$ as
$$A(G)=\left[\begin{array}{ccccc}
0&\alpha_1^T&\alpha_2^T&\cdots&\alpha_k^T\\
\alpha_1&A(G_1)&0&\cdots&0\\
\alpha_2&0&A(G_2)&\cdots&0\\
\vdots&\vdots&\vdots&\ddots&\vdots\\
\alpha_k&0&0&\cdots&A(G_k)
\end{array}\right],$$
where the left upper $0$ corresponds to the cut vertex $x$.
Observe that for each $i$ the equaiton $A(G_i)X=\alpha_i$ has a solution $\xi_i$; otherwise $r(G_i+x)=r(G_i)+2$; a contradiction.
Taking $Q$ as the following matrix with the same partition as $A(G)$,
$$Q=\left[\begin{array}{ccccc}
1&0&0&\cdots&0\\
-\xi_1&I&0&\cdots&0\\
-\xi_2&0&I&\cdots&0\\
\vdots&\vdots&\vdots&\ddots&\vdots\\
-\xi_k&0&0&\cdots&I
\end{array}\right],$$
we have
$Q^TA(G)Q=a \oplus A(G_1)\oplus A(G_2)\oplus \cdots \oplus A(G_k)$, where $a=-\sum_{i=1}^k\alpha_i\xi_i$.
The assumption $s(G) = s(G-x)+1$ implies that $a > 0$.
In particular, their exists some $l$ such that $\alpha_l\xi_l < 0$.
So, $A(G_l+x)$ is congruent to $(-\alpha_l\xi_l) \oplus A(G_1))$, which implies $s(G_l+x)=s(G_l)+1$.
Therefore $s(G)=s(G_l+x)+\sum_{j\neq l}s(G_{j})$. \hfill $\blacksquare$

\begin{coro}\label{cut vertex c5}
Let $x$ be a cut vertex of a graph $G $ and let $G_1, G_2, \ldots, G_k$  be all components of $G -x$.
If $s(G_i)\leq c_5(G_i)$ and $s(G_i+x)\leq c_5(G_i+x)$  for all $i$'s, then $s(G)\leq c_5(G)$.
\end{coro}

\noindent{\bf Proof.}
By Lemma \ref{vertex-deleting sign 2}, $s(G)\leq s(G-x)+1$.
If $s(G)\leq s(G-x)$, noting that $s(G-x) =\sum_{i=1}^k s(G_i)$ and $s(G_i)\leq c_5(G_i)$ for all $i$'s,  so we have $s(G)\leq\sum_{i=1}^k c_5(G_i)\leq c_5(G)$.
If $s(G)=s(G-x)+1$, by Lemma \ref{cut vertex sign 3}, $s(G)=s(G_l+x)+\sum_{j\ne l} s(G_j)$ for some $l$.
By the assumption for each $G_i$ and $G_i+x$,  we have $s(G)\leq c_5(G_l+x)+ \sum_{j\ne l} c_5(G_j)\leq c_5(G)$.
\hfill $\blacksquare$

\section{Signature of line graphs}
The {\it line graph} of a graph $G$, denoted by $L_G$, is the graph whose vertex set is $E(G)$, where two
vertices of $L_G$ are adjacent if and only if the corresponding edges are incident in $G$.

\begin{lemma} \label{s-graph}
If $G$ is one of the following graphs: a cycle with two pendant edges,
two cycles sharing a common vertex,  two cycles sharing a common path of length at least $1$,
where all cycles have length $2$ modulo $4$,
then $s(L_G) \le c_5(L_G)$.
\end{lemma}

\noindent{\bf Proof.}
First suppose $G$ is a cycle $C$ of length $2$ modulo $4$ with two pendant edges $e_1=x_1y_1$ and $e_2=x_2y_2$, where $y_1,y_2$ are pendant vertices of $G$.
If $x_1=x_2$ or $x_1, x_2$ are connected by paths on $C$ of even length, by Lemma \ref{luni}, $L_G$ is nonsingular.
Note that $L_G$ has an even order so that $s(L_G)$ is an even number.
By Theorem \ref{tree}, $s(L_{C}+e_1) \le 0$ as $L_{C}+e_1$ is bicyclic.
Now by Lemma \ref{vertex-deleting sign 2}, $s(L_G) \le s(L_{C}+e_1)+1 \le 1$.
So, $s(L_G) \le 0 =c_5(L_G)$.

If $x_1,x_2$ are  connected by paths on $C$ of odd length, then $\eta(L_G)=1$ by Lemma \ref{luni}.
Note that $C-x_1-x_2$ consists of two disjoint paths $P_1,P_2$ both with order $0$ or $2$ modulo $4$.
By Lemma \ref{sper}, it suffices to consider the line graphs $G_1,G_2$ in Fig. 3.1.
We have $s(G_1)=s(G_2)=-1$ by using {\scshape Mathematica}.

Next we consider the case that $G$ is two cycles sharing a common vertex.
Also by Lemma \ref{sper} it suffices to consider the line graph $G_3$ in Fig. 3.1.
By a direct calculation, we have $s(G_3)=-1$.

Finally we consider the case that $G$ is two cycles sharing a common path $P$ of length at least $1$.
We stress all cycles have length $2$ modulo $4$.
 If the path $P$ has length $1$, then by Lemma \ref{sper} it suffices to consider the line graph $G_4$ in Fig. 3.1.
By a direct calculation, we have $s(G_4)=-1$.
If $P$ has length greater than $1$, then by Lemma \ref{sper} it suffices to consider the line graphs $G_5,G_6$ in Fig. 3.1.
Also by calculation, we get $s(G_5)=s(G_6)=-1$. \hfill $\blacksquare$

\begin{center}
\vspace{3mm}
\includegraphics[scale=.6]{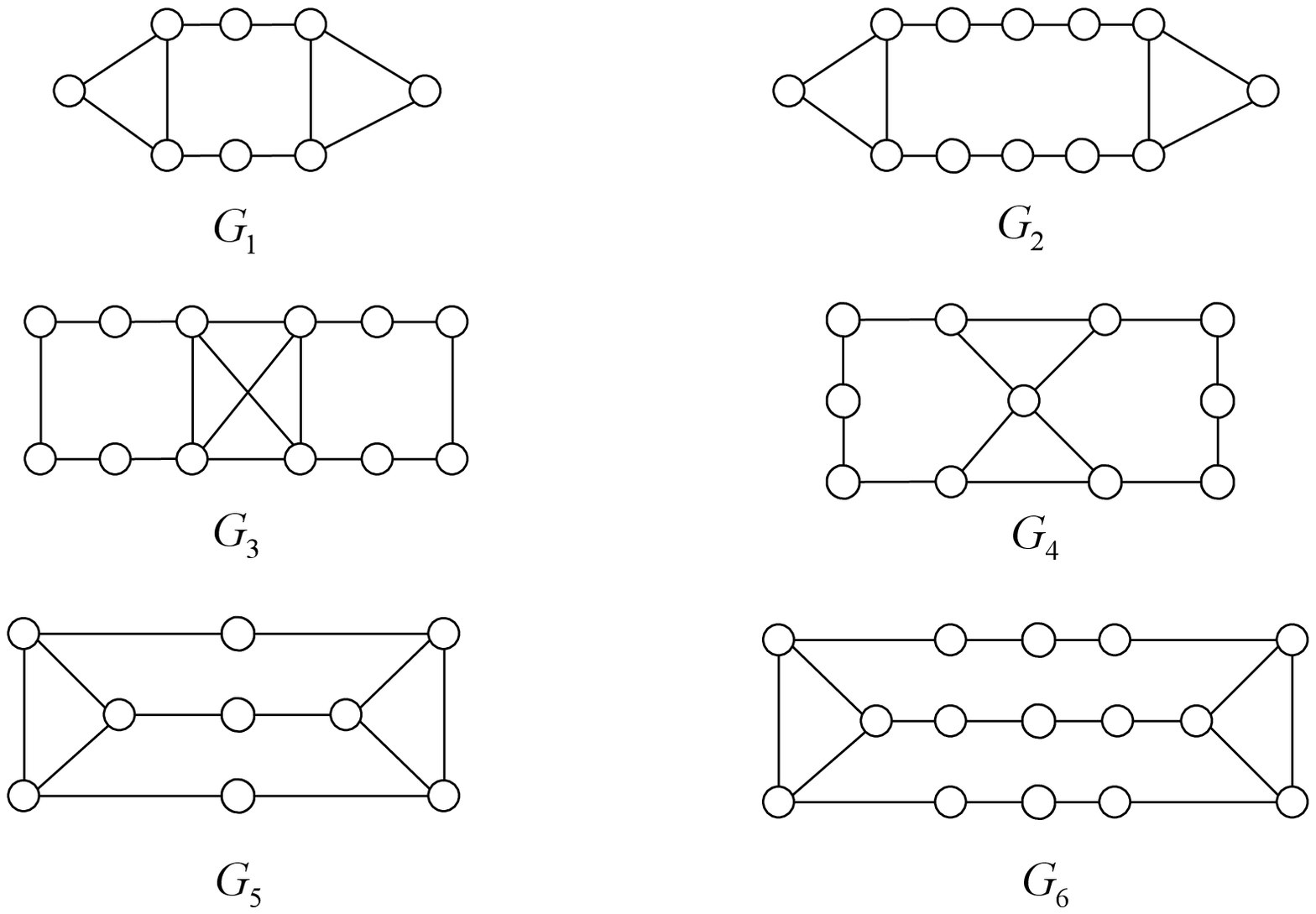}

{\small Fig. 3.1. The graphs in the proof of Lemma \ref{s-graph}}
\end{center}

\begin{theorem}\label{tree c5} Let $T$ be a tree with at least one edge, then $s(L_T)\leq c_5(L_T)$.
\end{theorem}

\noindent{\bf Proof.}
We use induction on the number of internal edges (i.e. non-pendant edges) of $T$ to prove the result.
If $T$ has no internal edges, then $T=K_{1,m}$ (i.e. a star), and hence $L_T=K_m$, a complete graph.
The result holds in this case by a simple verification.

Suppose the result holds for all trees with $k \;(\geq 0)$ internal edges.
Let $T$ be a tree with $k+1$ internal edges and let $e$ be one of the internal edges of $T$.
Then $T-e$ consists of two subtrees $T_1,T_2$ of $T$.
Obviously, each  $T_i$  and each $T_i+e$ has fewer internal edges than that of $T$.
By induction we have $s(L_{T_i})\leq c_5(L_{T_i})$ and $s(L_{T_i}+e)\leq c_5(L_{T_i}+e)$ for each $i=1,2$.
Noting that $e$ is a cut vertex of  $L_T$, so $s(L_T)\leq c_5(L_T)$ by Corollary \ref{cut vertex c5}. \hfill $\blacksquare$

\begin{theorem}\label{theorem 1}
Let $G$ be a graph without isolated vertices. Then $-c_3(L_G)\leq s(L_G)\leq c_5(L_G)$.
\end{theorem}

\noindent{\bf Proof.}
Without loss of generality we may assume $G$ is connected.
Let $\Theta(G)$ be the set of edges of $G$ with at least one endpoint having degree greater than $2$, and let $\theta(G):=|\Theta(G)|$.
We will use induction on $\theta(G)$ to prove the left inequality.
If $\theta(G)=0$, namely each vertex of $G$ has degree $1$ or $2$, then $G$ is the disjoint union of paths and/or cycles.
Thus, $L_G$ is the disjoint union of paths and/or cycles. By Theorem \ref{tree}, we have $-c_3(L_G)\leq s(L_G)$.

Assume that $-c_3(L_H)\leq s(L_H)$ for all graphs $H$ with $\theta(H)\leq k$, where $ k \ge 0$.
Let $G$ be a graph with $\theta(G)=k+1$ and let $x$ be a vertex of $G$ with degree at least $3$.
Suppose $e$ is an edge incident to $x$. Then the vertex $e$ of $L_G$ is contained in one triangle.
So $c_3(L_{G-e})=c_3(L_G-e) \leq c_3(L_G)-1$.
By Lemma \ref{vertex-deleting sign 2} and by induction,
$$ s(L_G)\geq s(L_G-e)-1 =s(L_{G-e})-1\geq -c_3(L_{G-e})-1\geq -c_3(L_G).$$

Next, set $o(G):=|E(G|-|V(G)|+1$, the {\it dimension} of $G$.
 We also use induction on $o(G)$ to prove the right inequality.
 If $o(G)=0$, then $G$ is a tree, and the result holds in this case by Theorem \ref{tree c5}.
 Assume the result holds for all connected graphs $G$ with $o(G)\leq k$, where $k \ge 0$.
 Let $G$ be a connected graph with $o(G)=k+1$.
Note that $G$ must contain cycles.
A cycle $C$ of $G$ is said of {\it type $l$} if there are exactly $l$ edges between $C$ and $G-C$.

{\bf Case 1:}  If $G$ contains a cycle $C$ of type $l$ with $l \ge 3$, letting $m$ be the length of $C$ and letting $e_1,e_2,e_3$ be three edges joining $C$ and $G-C$, then the line graphs $L_C, L_C+e_1, L_C+e_1+e_2,L_C+e_1+e_2+e_3$ contain cycles of length $m,m+1,m+2,m+3$ respectively.
Surely one cycle among them must have length $1$ modulo $4$. 
Deleting an arbitrary edge, say $e$ on the cycle $C$, will break the cycle of length $1$ modulo $4$ and decrease the dimension of $G$.
That is, $c_5(L_G-e)\leq c_5(L_G)-1$, and $o(G-e) < o(G)$.
Now by Lemma \ref{vertex-deleting sign 2} and by induction,
$$s(L_G)\leq s(L_G-e)+1 =s(L_{G-e})+1\leq c_5(L_{G-e})+1 \leq c_5(L_G).$$

{\bf Case 2:} If $G$ contains a cycle of type $1$, say $C$, then $C$ is connected to $G-C$ by an edge, say $e=xy$, where $x \in V(C)$ and
$y \in V(G-C)$.
Surely $e$ is a cut edge of $G$.
If $G=C+y$, then $L_G$ is bicyclic and the result holds by Theorem \ref{tree}.
If $G \ne C+y$, then $e$ is a cut vertex of $L_G$, $G-e$ has two components: $C$ and another subgraph say $D$, where $o(D) < o(G)$ and $o(D+x) < o(G)$.
So, by induction, $s(L_D) \le c_5(L_D)$ and $s(L_D+e) \le c_5(L_D+e)$.
Observe that $s(L_C) \le c_5(L_C)$ and $s(L_C+e) \le c_5(L_C+e)$ by Theorem \ref{tree}.
The result now follows by Corollary \ref{cut vertex c5}.

{\bf Case 3:} If all cycles of $G$ are of type $2$, then $G$ is either (i) one obtained from a cycle with two pendant edges (denoted by $H$) by possibly attaching trees at the pendant vertices of $H$, or (ii) two cycle sharing a common vertex or a common path of length at least $2$,
or (iii) $G$ is obtained from a tree by replacing some vertices of degree $2$ by cycles.

If $G$ is one of graphs in (i) and (ii), and in addition if one cycle has odd length or length $0$ modulo $4$, then we will find a cycle in $G$ of length $1$ modulo $4$ containing
the edges of the cycle.
Similar to Case 1, deleting an arbitrary edge on the cycle will break the cycle of length $1$ modulo $4$ and decrease the dimension of $G$.
The result will follows by Lemma \ref{vertex-deleting sign 2} and by induction.

Now assume $G$ is one of graphs in (i) and (ii), and  all cycles have length $2$ modulo $4$.
If $G$ is exactly the graph $H$ (a special case of (i)) or a graph in (ii), we get the result by Lemma \ref{s-graph}.

If $G$ is a graph in (i) obtained from $H$ by attaching exactly one tree $T$ at the pendant vertex of a pendant edge say $e$, then
$G$ contains a cut edge say $e$ such that $G-e$ has two components: $G_1,T$, where $G_1$ is the cycle together with a pendant edge.
Note that $e$ is a cut vertex of $L_G$, and $s(L_{T}) \le c_5(L_{T})$, $s(L_{T}+e) \le c_5(L_{T}+e)$ by Theorem \ref{tree c5}.
Also $s(L_{G_1}) \le c_5(L_{G_1})$ by Theorem \ref{tree} as $L_{G_1}$ is bicyclic, $s(L_{G_1}+e) \le c_5(L_{G_1}+e)$ by Lemma \ref{s-graph}.
The result now follows by Corollary \ref{cut vertex c5}.

If $G$ is a graph in (i) obtained from $H$ by attaching two trees at the pendant vertices of two pendant edge say $e_1,e_2$ respectively, 
Then $G-e_2$ has two components: $G_1,G_2$, where $G_1$ contains the cycle and $G_2$ is a tree.
Note that in the graph $G_1$ the cycle is of type $1$, and hence $s(L_{G_1}) \le c_5(L_{G_1})$ by the result in Case 2.
Also $s(L_{G_1}+e) \le c_5(L_{G_1}+e)$ by what we have proved in this case.
So the result also follows by Corollary \ref{cut vertex c5}.

If $G$ is a graph in (iii), then there exists a cut edge $e$ of $G$ such that $G-e$ has two components: $G_1,G_2$, where $G_1,G_2$ both contain cycles.
So $o(G_i)<o(G), o(G_i+e)< o(G)$, and by induction $s(L_{G_i})<c_5(L_{G_i}), s(L_{G_i+e})< c_5(L_{G_i+e})$ for $i=1,2$.
Note that $e$ is a  cut vertex of $L_G$.
The result also follows by Corollary \ref{cut vertex c5}.

{\bf Case 4:} If $G$ contains a cycle of type $0$ and contains no cycles of type $1$ or type $l$ with $l \ge 3$,  then $G$ itself is the cycle or the cycle with a chord (an edge with two endpoints on the cycle).
Clearly the result holds if $G$ is a cycle.
If $G$ is a cycle with a chord, letting $C_1,C_2$ be two smaller cycles containing the chord, if one cycle has odd length or length $0$ modulo $4$,
then the result follows by a similar discussion as in Case 3.
Otherwise, $C_1,C_2$, and hence $C$ all have length $2$ modulo $4$.
In this case, we can get the result by Lemma \ref{s-graph}. \hfill $\blacksquare$

\section{Signature of power trees}
Recall that the {\it $k$-th power} $G^{k}$ of a graph $G$ is obtained from $G$ by adding edges between all pairs of vertices within distance at most $k$.
In particular $G^1$ is exactly the graph $G$, and $G^2$ is called the {\it square} of $G$.

%

\begin{lemma} \label{power35}Let $G$ be a graph on at least $5$ vertices. If $k\geq 2$, then in the graph $G^{k}$ every vertex $v$ is contained in at least one $C_3$ and one $C_5$. That is to say, $c_{3}(G^{k}-v)\leq c_{3}(G^{k})-1$ and $c_{5}(G^{k}-v)\leq c_{5}(G^{k})-1$.
\end{lemma}

\noindent{\bf Proof.} Let $H$ be an arbitrary connected graph induced by five vertices of $G$. Then $H$ contains one of $H_{1},H_{2},H_{3}$ as a subgraph; see Fig. 4.1.
Thus $G^{2}$, and hence $G^k$ contains $H_{1}^{2}$ as a subgraph by considering the squares of $H_{1},H_{2},H_{3}$.
Note that in $H_{1}^{2}$ each vertex is contained in at least one  $C_3$ and one $C_5$.
The result follows.\hfill $\blacksquare$


\begin{center}
\vspace{3mm}
\includegraphics[scale=.6]{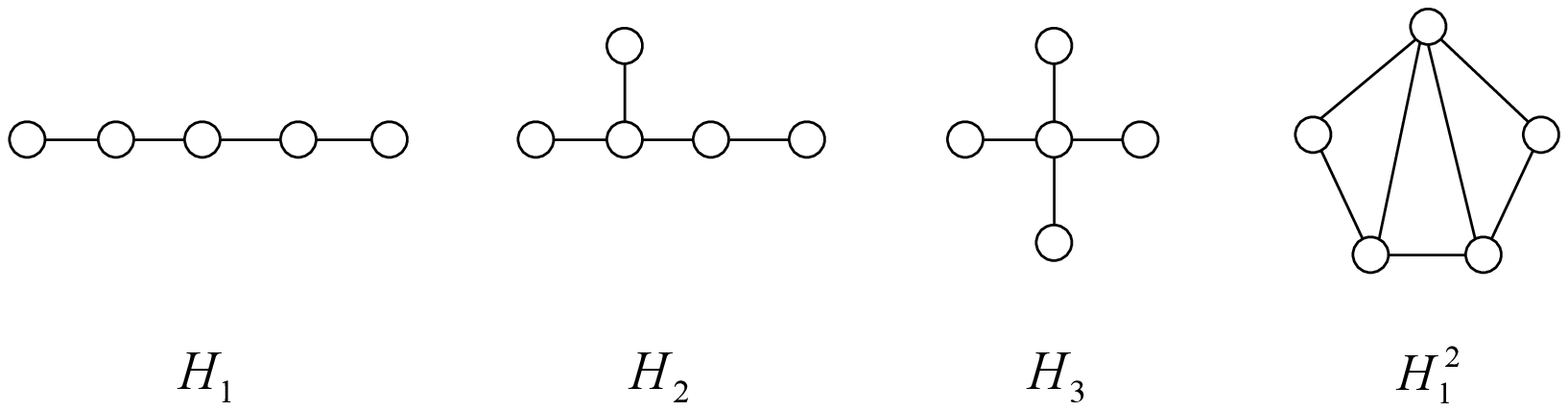}

{\small Fig. 4.1. The graphs in the proof of Lemma \ref{power35}}
\end{center}

\begin{theorem}\label{theorem 2}
If $G$ is a tree, then $-c_3(G^{k})\leq s(G^{k})\leq c_5(G^k)$ for $k\geq 2$.
\end{theorem}

\noindent{\bf Proof.}
If $|V(G)|\leq 4$, the result follows by a direct calculation.
Assume the result holds for all trees on $n$ vertices, where $n \ge 4$.
Let $G$ be a tree on $n+1$ vertices.
By Lemma \ref{power35}, $c_{3}(G^k-v)\leq c_{3}(G^k)-1$ and $c_{5}(G^k-v)\leq c_{5}(G^k)-1$ for an arbitrary vertex $v$ of $G$.
Let $u$ be a pendant vertex of $G$. Then $G^k-u=(G-u)^k$.
So Lemma \ref{vertex-deleting sign 2} and by induction,
$$s(G^k) \le s(G^k-u)+1 =s((G-u)^k)+1 \le c_5((G-u)^k)+1=c_5(G^k-u)+1 \le c_5(G^k)-1+1=c_5(G^k).$$
Similarly,
$$ s(G^k) \ge s((G-u)^k)-1 \ge -(c_3(G^k-u)+1) \ge -c_3(G^k).$$
The result follows.\hfill $\blacksquare$

%
%
%
%
%
%

Recall that the {\it total graph} $T_G$ of $G$ is the graph with vertex set corresponding to union of vertex and edge sets of $G$, with two
vertices of $T_G$ adjacent if and only if the corresponding elements in $G$ are adjacent or incident.
It is known that $T_G = S(G)^{2}$ (or see \cite{har}), where $S(G)$ is the subdivision of $G$.
If $G$ is a tree, then $S(G)$ is also a tree.
So we have the following corollary.

\begin{coro} If $G$ is a tree, then $-c_{3}(T_G)\leq s(T_G)\leq c_{5}(T_G)$.
\end{coro}


{\small
\begin{thebibliography}{90}


\bibitem{atk} P. W. Atkins, J. de Paula, {\it Physical Chemistry}, eighth ed., Oxford University Press, 2006.


\bibitem{col} L. Collatz, U. Sinogowitz, Spektren endlicher Grafen, {\it Abh. Math. Sem. Univ. Hamburg,}  21 (1957) 63-77.

\bibitem{cve} D. Cvetkov\'{\i}c, M. Doob, H. Sachs, {\it Spectra of Graphs - Theory and Application}, Academic Press, New York, 1980.

\bibitem{gon} S. C. Gong, G. H. Xu, On the nullity of a graph with cut-points, {\it Linear Algebra Appl.,} 436 (2012) 135-142.

\bibitem{gut} I. Gutman, I. Sciriha, On the nullity of line graphs of trees, {\it Discrete Math.,} 232 (2001) 35-45.

\bibitem{har} F. Harary, {\it Graph Theory}, Addison-Wesley, Reading, 1969.

\bibitem{li} H.-H. Li, Y.-Z. Fan,  L. Su,  On the nullity of the line graph of unicyclic graph with depth one, {\it Linear Algebra Appl.,}  437 (2012) 2038-2055.

\bibitem{lisc} S. Li, F. Song, On the positive and negative inertia of weighted graphs, {\it arXiv}: 1307.5110.

\bibitem{ma} H. Ma, W. Yang, S. Li, Positive and negative inertia index of a graph, {\it Linear Algebra Appl.,} 438 (2013) 331-341.

\bibitem{maxb} X. Ma, D. Wong, M. Zhu, The positive and the negative inertia index of line graphs of trees,  {\it Linear Algebra
Appl.} (2013), http://dx.doi.org/10.1016/j.laa.2013.08.024.


\bibitem{sci} I. Sciriha, On singular line graphs of trees, {\it Congr. Numer.,} 135 (1998) 73-91.


\bibitem{yu} G. Yu,  L. Feng, Q. Wang, Bicyclic graphs with small positive index of inertia,  {\it Linear Algebra Appl.,} 438 (2013) 2036-2045.

\bibitem{yugh} G. Yu, X.-D. Zhang, L. Feng, The inertia of weighted unicyclic graphs, {\it arXiv}: 1307.0059.
\end{thebibliography}
}
\end{document}